\documentclass[10pt]{article}
\usepackage[letterpaper]{geometry}
\usepackage{hicss}
\usepackage{times}
\makeatletter
\long\def\@caption#1[#2]#3{%
  \par
  \addcontentsline{\csname ext@#1\endcsname}{#1}{%
    \protect\numberline{\csname the#1\endcsname}{\ignorespaces #2}}%
  \begingroup
    \@parboxrestore
    \normalsize
    \@makecaption{\csname fnum@#1\endcsname}{\ignorespaces #3}\par
  \endgroup
}
\long\def\@makecaption#1#2{%
  \vskip 4pt%
  \setbox\@tempboxa\hbox{ \fontfamily{cmss}\selectfont\textbf{\small #1.~ }}%
  \setlength\captionindent{\wd\@tempboxa}\divide\captionindent by 2%
  \parbox[t]{\hsize}{\centering\fontfamily{cmss}\selectfont\textbf{\captionsize
    \hangindent\captionindent \unhbox\@tempboxa#2}}%
}
\makeatother
\usepackage[none]{hyphenat}
\usepackage{url}
\usepackage{latexsym}
\usepackage{indentfirst}
\usepackage{amsmath,amssymb,amsfonts}
\usepackage{bm}
\usepackage{booktabs}
\usepackage{multirow}
\usepackage{array}
\usepackage{graphicx}
\graphicspath{{images/}}
\usepackage{algorithm}
\usepackage[noend]{algpseudocode}
\usepackage[hidelinks]{hyperref}
\usepackage[
  backend=biber,
  style=numeric-comp,
  sorting=none,
  giveninits=true,
  hyperref=true,
]{biblatex}
\usepackage{tikz}
\usepackage{pgfplots}
\pgfplotsset{compat=1.18}
\usepackage{xcolor}
\newcommand{\bv}[1]{\mathbf{#1}}   
\newcommand{\vect}[1]{\bm{#1}}     
\newcommand{\lb}[1]{\underline{#1}}
\newcommand{\ub}[1]{\overline{#1}}

\newcommand{\im}{\mathbf{j}}  

\newcommand{\NODES}{\mathcal{N}}
\newcommand{\EDGES}{\mathcal{E}}

\newcommand{\GENERATORS}{\mathcal{G}}

\newcommand{\CBLOC}{\NODES^{\mathrm{CB}}}
\newcommand{\TAPLOC}{\EDGES^{\mathrm{T}}}
\newcommand{\CBPOSSET}{\mathcal{C}}
\newcommand{\TAPPOSSET}{\mathcal{T}_{\mathrm{pos}}}
\newcommand{\DEVSET}{\mathcal{D}}

\newcommand{\Eiplus}{\EDGES_i^{+}}
\newcommand{\Eiminus}{\EDGES_i^{-}}

\newcommand{\V}{\bv{V}}
\newcommand{\VM}{\bv{v}}
\newcommand{\VA}{\vect{\theta}}

\newcommand{\SG}{\bv{S}^{\mathrm{g}}}
\newcommand{\PG}{\bv{p}^{\mathrm{g}}}
\newcommand{\QG}{\bv{q}^{\mathrm{g}}}
\newcommand{\SF}{\bv{S}^{\mathrm{f}}}
\newcommand{\ST}{\bv{S}^{\mathrm{t}}}
\newcommand{\PF}{\bv{p}^{\mathrm{f}}}
\newcommand{\QF}{\bv{q}^{\mathrm{f}}}
\newcommand{\PT}{\bv{p}^{\mathrm{t}}}
\newcommand{\QT}{\bv{q}^{\mathrm{t}}}
\newcommand{\DEV}{\bv{d}}
\newcommand{\CDEV}{\bv{c}}
\newcommand{\TDEV}{\bv{t}}

\newcommand{\CB}{\beta}
\newcommand{\TAP}{\tau}
\newcommand{\Ybr}{\mathrm{Y}^{\mathrm{br}}}

\newcommand{\XACOPF}{\bv{x}_{\mathrm{acopf}}}
\newcommand{\XACPF}{\bv{x}_{\mathrm{acpf}}}
\newcommand{\XVVO}{\bv{x}_{\mathrm{vvo}}}
\newcommand{\XVVOACOPF}{\bv{x}_{\mathrm{vvo|acopf}}}
\newcommand{\XVVOACPF}{\bv{x}_{\mathrm{vvo|acpf}}}

\newcommand{\vm}{v}

\newcommand{\bij}{b_{ij}}
\newcommand{\yff}{y^{\mathrm{ff}}}
\newcommand{\yft}{y^{\mathrm{ft}}}
\newcommand{\ytf}{y^{\mathrm{tf}}}
\newcommand{\ytt}{y^{\mathrm{tt}}}
\newcommand{\smax}{\overline{s}}
\newcommand{\gs}{g^{\mathrm{s}}}
\newcommand{\bs}{b^{\mathrm{s}}}
\newcommand{\ys}{y^{\mathrm{s}}}
\newcommand{\Sd}{S^{\mathrm{d}}}
\newcommand{\pd}{p^{\mathrm{d}}}
\newcommand{\qd}{q^{\mathrm{d}}}
\newcommand{\cost}{c}
\newcommand{\bcb}{b^{\mathrm{cb}}}

\newcommand{\vref}{v^{\mathrm{ref}}}
\newcommand{\qgref}{q^{\mathrm{ref}}}
\newcommand{\pgref}{p^{\mathrm{ref}}}

\newcommand{\pgmin}{\lb{p}^{\mathrm{g}}}
\newcommand{\pgmax}{\ub{p}^{\mathrm{g}}}
\newcommand{\qgmin}{\lb{q}^{\mathrm{g}}}
\newcommand{\qgmax}{\ub{q}^{\mathrm{g}}}
\newcommand{\vmmin}{\lb{v}}
\newcommand{\vmmax}{\ub{v}}
\newcommand{\dvamin}{\lb{\Delta}\theta}
\newcommand{\dvamax}{\ub{\Delta}\theta}

\newcommand{\lamv}{\lambda_v}
\newcommand{\lamq}{\lambda_q}
\newcommand{\lamp}{\lambda_p}
\newcommand{\lamc}{\lambda_{c}}

\newcommand{\thetaref}{\VA_{\mathrm{ref}}}

\newcounter{model}
\newcommand{\modeltitle}[2]{%
  \par\noindent\rule{\columnwidth}{0.4pt}\\[2pt]%
  \noindent\refstepcounter{model}\label{#1}\textbf{Model~\themodel} \hspace{0.3em}(#2)\\[-0.25em]%
  \noindent\makebox[\columnwidth][l]{\rule{\columnwidth}{0.4pt}}\\[-1.1em]%
}
\newcommand{\modelbottomrule}{%
  \vspace{-1.25em}%
  \noindent\rule{\columnwidth}{0.4pt}\par\vspace{0.25em}%
}
\newcommand{\modeltightbottomrule}{%
  \vspace{-2.00em}%
  \noindent\rule{\columnwidth}{0.4pt}\par\vspace{2em}%
}
\newcommand{\modelhomotopybottomrule}{%
  \vspace{-30pt}%
  \noindent\rule{\columnwidth}{0.4pt}\par\vspace{39pt}%
}

\title{Improving Stability and Economic Operation in Transmission Systems through Volt/VAR Optimization}
\author{Shuaicheng Tong \\
 School of Industrial and Systems Engineering \\
 Georgia Institute of Technology \\
 {\underline{ stong38@gatech.edu}}  \\ \And
 Pascal Van Hentenryck \\
 School of Industrial and Systems Engineering \\
 Georgia Institute of Technology \\
 {\underline{ pvh@gatech.edu} } \\ }
\date{}

\begin{document}
\maketitle

\begin{abstract}
Transmission system operators often reconcile market-cleared DC dispatches with
AC physics through power-flow solves (ACPF), yet the resulting setpoints
can still violate voltage (Volt) and reactive-power (VAR) limits.
Maintaining secure voltage profiles and adequate VAR support therefore requires
fast corrective decisions that are implementable in operation.
This paper presents a novel {\em homotopy-based continuation method} for discrete-control Volt/VAR Optimization (VVO) that coordinates switchable devices, such as on-load tap-changing transformers (OLTCs) and capacitor banks (CBs).
Experiments on IEEE, PEGASE, and RTE systems show that the proposed VVO produces
AC-feasible setpoints within practical runtime, while reducing voltage
deviation, VAR dispatch, and generation cost.
VVO also achieves comparable performance using ACPF-adjusted DC dispatches as
inputs relative to AC-feasible dispatches, indicating that it can be integrated
naturally into existing transmission dispatch practices to strengthen grid
operation.
\end{abstract}

\section{Introduction}
Transmission system operators (TSOs) must deliver active power economically
while keeping voltage, reactive power, and branch flows within secure limits.
Voltage (Volt) and reactive-power (VAR) control is challenging because voltage
deviations and insufficient VAR support reduce operational flexibility and
increase stability risks, especially under stressed system conditions
\cite{Kundur2004_StabilityDefs,ElizondoSurvey2017_TransVoltageControl}.

However, upstream market clearing and generator dispatch commonly rely on DC
optimal power flow (DCOPF) models, which linearize the nonconvex constraints in
AC optimal power flow (ACOPF) by omitting voltage magnitude and reactive power
variables and constraints. This creates a modeling gap between the dispatch setpoints
TSOs receive and the AC physics they must ultimately regulate. As a result,
DCOPF setpoints are not AC-feasible by themselves. Indeed, Baker showed that,
under mild assumptions, the feasible regions of DCOPF and ACOPF are disjoint
\cite{baker2021solutions}.

A common feasibility restoration technique is to run an AC power flow (ACPF) on
DCOPF setpoints. This recovers AC variables, but the resulting solution may still 
violate voltage, reactive power, and thermal limits. Recent work by Boateng et al. explores DCOPF and
ACPF variants that substantially reduce violations through distributed slack
and PV/PQ bus type switching~\cite{boateng2025towards}. The limitation remains
fundamental to ACPF-based recovery: it does not guarantee an AC-feasible setpoint or
optimize operational objectives. It can also become computationally expensive on large systems,
with reported runtimes exceeding ACOPF runtimes.

These observations motivate a scalable optimization layer that coordinates
available controls and explicitly enforces AC feasibility. In transmission
networks, this leads to a Volt/VAR Optimization (VVO) problem over control
devices such as on-load tap-changing transformers (OLTCs) and capacitor 
banks (CBs), which adjust voltage magnitudes and reactive power
flows, respectively \cite{EPRI2024VVO}. Coordinating these devices is especially relevant 
because their settings provide additional degrees of freedom for feasibility
recovery and setpoint improvement.

Optimizing the discrete settings of OLTCs and CBs under nonconvex AC
power flow constraints results in a mixed-integer nonlinear program (MINLP),
making the VVO problem difficult to solve directly at transmission scale
\cite{PapalexopoulosLargeOPF, CAPITANESCU201657}. Prior work has
proposed tractable methods for incorporating discrete controls into
reactive power optimization. Classical approaches
incorporated discrete shunts and transformer controls into Newton or
interior-point methods \cite{Liu1992,Acha2000,Liu2002ExtendedIP}. Sensitivity-
and rounding-based methods use local information from continuous OPF solutions
to rank discrete adjustments
\cite{Capitanescu2010SensitivityOPF,Macfie2010ShuntRounding}. Homotopy methods
solve a sequence of relaxed continuous problems that gradually drive discrete
controls toward admissible settings \cite{McNamara2022TwoStageHomotopy}. Other
approaches include semidefinite-relaxation branch-and-bound methods for
reactive-power dispatch with discrete controllers and global mixed-integer
ACOPF algorithms based on adaptive relaxations and cutting planes
\cite{Constante2021_TPWRS_ORPD_SDP_BB,Aigner2023ACOPFGlobalOptimality}.

A recent paper formulates the VVO problem as a MINLP that co-optimizes OLTC and
CB settings, solving it through a relax--round--resolve heuristic that scales
to large transmission test cases \cite{tong2026vvo}. In that study, however,
feasibility had already been enforced upstream: the input setpoints were ACOPF
solutions, and VVO was used to improve the setpoints. This paper considers a
more realistic operation setting in which the input operating point comes from a
DCOPF$\rightarrow$ACPF pipeline that is not AC-feasible. In this setting, VVO acts
as a post-dispatch corrective optimization layer that optimizes Volt/VAR
quality and generation cost while enforcing AC feasibility.

This paper makes three main contributions. First, it improves the realism of 
a transmission VVO model by limiting excessive device movements. Second, it
strengthens numerical stability by embedding a homotopy continuation stage
that moves gradually from relaxed device positions to implementable settings.
Third, it evaluates VVO under a more realistic operational setting in which the
input comes from a DCOPF $\rightarrow$ ACPF setpoint, demonstrating
the proposed approach's scalability and effectiveness in feasibility recovery,
Volt/VAR improvement, and cost reduction.

\section{Problem Formulation}\label{sec:model}

\subsection{Network Notation}

Consider a power network represented as a directed graph
$G\,{=}\,(\NODES,\EDGES)$, where the set of buses and the set of branches
are denoted as $\NODES\,{=}\,\{1,\ldots,N\}$ and
$\EDGES\subseteq\NODES\times\NODES$, respectively.
Let $\Eiplus$ and $\Eiminus$ be the sets of branches leaving and entering bus $i$ and let $\GENERATORS$ be the set of generators.
For notational simplicity, the model is stated for one generator and one load per bus with no
parallel branches. The implementation and experimental workflows support parallel branches and arbitrary
numbers of loads and devices.

Let $\im$ be the imaginary unit and $(\cdot)^{\star}$ be the complex conjugation.
Within the model formulations, bold symbols denote decision variables, with the
exception of the imaginary unit $\im$.
At bus $i\in\NODES$, the complex voltage, power generation, and power demand are
$\V_{i}\,{=}\,\VM_{i}\angle\VA_{i}$,
$\SG_{i}\,{=}\,\PG_{i}\,{+}\,\im\QG_{i}$,
and $\Sd_{i}\,{=}\,\pd_{i}\,{+}\,\im\qd_{i}$, respectively.
On branch $ij\in\EDGES$, the complex power injections at the from- and to-end terminals are
$\SF_{ij}\,{=}\,\PF_{ij}\,{+}\,\im\QF_{ij}$
and $\ST_{ij}\,{=}\,\PT_{ij}\,{+}\,\im\QT_{ij}$, respectively.

\subsection{Volt/VAR Optimization Model}
VVO coordinates controls of OLTC tap positions and CB module offsets jointly with voltage,
generation, and branch flow variables under nonconvex AC power flow equations and
operating limits. Let \(\CDEV_i\) denote the CB module offset, where \(\CDEV_i=0\) leaves the original shunt
susceptance unchanged. Let \(\TDEV_{ij}\) denote the tap position, with reference position 
\(\TDEV^0_{ij}\) from the input case files. Let \(\CBPOSSET\subset\mathbb{Z}\) and \(\TAPPOSSET\subset\mathbb{Z}\)
denote the admissible integer sets for CB offsets and tap positions, respectively. Let \(\CBLOC\) and \(\TAPLOC\) denote the
set of buses equipped with CBs and the set of branches with transformers. The
device-dependent shunt and tap quantities entering the power flow equations are
\begin{align}
  \CB_i(\CDEV_i) &= \CDEV_i\bcb \in \mathbb{R}, \quad i\in\CBLOC, \label{eq:cb-map}\\
  \TAP_{ij}(\TDEV_{ij}) &= 1+\Delta\tau\,\TDEV_{ij} \in \mathbb{R}, \quad ij\in\TAPLOC, \label{eq:tap-map}
\end{align}
where \(\bcb\) is the per module susceptance value and \(\Delta\tau\)
is the tap step size. For non-CB buses, \(\CB_i(\CDEV_i)=0\). For branches without OLTCs, \(\TAP_{ij}(\TDEV_{ij})=1\).
To avoid excessive switching and associated maintenance costs, device movements
are limited by step budgets
\[
  B_{\beta}=\lceil\rho_{\beta}|\CBLOC|\rceil,
  \qquad
  B_{\tau}=\lceil\rho_{\tau}|\TAPLOC|\rceil,
\]
where \(\rho_{\beta}\) and \(\rho_{\tau}\) are the CB and tap budget multipliers.

The total shunt admittance at bus $i$ is
\begin{equation}
  \ys_i(\CDEV_i) = \gs_i + \im\bigl(\bs_i+\CB_i(\CDEV_i)\bigr),
  \label{eq:shunt}
\end{equation}
where $\gs_i$ and $\bs_i$ are the fixed shunt conductance and susceptance.
For compactness, let $a_{ij}=\TAP_{ij}(\TDEV_{ij})$ denote the realized tap
ratio. For a branch $ij$ with a transformer, the branch admittance submatrix is
\begin{equation}
  \Ybr_{ij}(a_{ij})=
  \begin{bmatrix}
    \yff_{ij} \times a_{ij}^{-2} & \yft_{ij} \times a_{ij}^{-1}\\
    \ytf_{ij} \times a_{ij}^{-1} & \ytt_{ij}
  \end{bmatrix}\in\mathbb{C}^{2\times2},
  \label{eq:Ybr}
\end{equation}
where $a_{ij}>0$, and phase shift angles are assumed to be fixed.
Kirchhoff's current law (KCL) enforcing nodal power balance at bus $i\in\NODES$ is
\begin{equation}
  \SG_i - \Sd_i - (\ys_i(\CDEV_i))^{\star}\VM_i^{2}
  = \sum_{e\in\Eiplus}\SF_e + \sum_{e\in\Eiminus}\ST_e.
  \label{eq:kcl}
\end{equation}
Ohm's law relating the branch terminal voltages to the forward and reverse flows:
\begin{align}
  \label{eq:ohm:complex}
  \begin{bmatrix}
    \SF_{ij}\\
    \ST_{ij}
  \end{bmatrix}
  &=
  \begin{bmatrix}
    \V_{i} & \\
    & \V_{j}
  \end{bmatrix}
  \times
  \Bigl(\Ybr_{ij}(a_{ij})\Bigr)^{\star}
  \times
  \begin{bmatrix}
    \V_{i}^{\star}\\
    \V_{j}^{\star}
  \end{bmatrix}
  ,
\end{align}
and thermal constraints read
\begin{align}
  \label{eq:thermal}
  |\SF_{ij}|,\;|\ST_{ij}| &\leq \smax_{ij},
  \qquad\forall\,ij\in\EDGES.
\end{align}

With these AC equations and device restrictions, the VVO objective is
\begin{equation}
  \begin{aligned}
  \psi(\VM,\QG,\PG)
    &=
    \lamv\psi_{v}(\VM)
    + \lamq\psi_{q}(\QG)\\
    &\quad
    + \lamp\psi_{p}(\PG)
    + \lamc\psi_{c}(\PG),
  \end{aligned}
  \label{eq:vvo_psi}
\end{equation}
where
\begin{align}
  \psi_{v}(\VM) &= \sum_{i\in\NODES}(\VM_{i}-\vref_{i})^{2},
  \label{eq:psiv}\\
  \psi_{q}(\QG) &= \sum_{i\in\GENERATORS}(\QG_{i}-\qgref_{i})^{2},
  \label{eq:psiq}\\
  \psi_{p}(\PG) &= \sum_{i\in\GENERATORS}(\PG_{i}-\pgref_{i})^{2},
  \label{eq:psip}\\
  \psi_{c}(\PG) &= \sum_{i\in\GENERATORS}\cost_{i}(\PG_{i}).\label{eq:psic}
\end{align}
Terms~\eqref{eq:psiv}--\eqref{eq:psip} penalize deviations from the voltage
target, reactive power reference, and upstream dispatch \(\pgref\).
\eqref{eq:psic} accounts for generation cost. The weights
\(\lamv,\lamq,\lamp,\lamc\) tune their relative importance.

\modeltitle{model:vvo}{Volt/VAR Optimization}
\begingroup
\postdisplaypenalty=10000
\begin{align}
  \min\quad
  &\psi(\VM,\QG,\PG) \label{eq:vvo_obj}\\
  \text{s.t.}\quad
  &\eqref{eq:kcl},\;\eqref{eq:ohm:complex},\;\eqref{eq:thermal},\nonumber\\
  &\thetaref=0, \label{eq:vvo_ref}\\
  &\pgmin_{i}\leq\PG_{i}\leq\pgmax_{i}, \quad\forall\,i\in\GENERATORS,
    \label{eq:vvo_pg}\\
  &\qgmin_{i}\leq\QG_{i}\leq\qgmax_{i}, \quad\forall\,i\in\GENERATORS,
    \label{eq:vvo_qg}\\
  &\dvamin_{ij}\leq\VA_{i}-\VA_{j}\leq\dvamax_{ij},
    \quad\forall\,ij\in\EDGES, \label{eq:vvo_ang}\\
  &\vmmin_{i}\leq\VM_{i}\leq\vmmax_{i}, \quad\forall\,i\in\NODES,
    \label{eq:vvo_vm}\\
  &\CDEV_{i}\in\CBPOSSET, \quad\forall\,i\in\CBLOC,
    \label{eq:vvo_cb}\\
  &\TDEV_{ij}\in\TAPPOSSET, \quad\forall\,ij\in\TAPLOC,
    \label{eq:vvo_tap}\\
  &\sum_{i\in\CBLOC}|\CDEV_i|\leq B_{\beta},
    \label{eq:vvo_cb_budget}\\
  &\sum_{ij\in\TAPLOC}|\TDEV_{ij}-\TDEV^0_{ij}|\leq B_{\tau}.
    \label{eq:vvo_tap_budget}
\end{align}
\modeltightbottomrule
\endgroup

The proposed VVO model is formulated as a MINLP in Model~\ref{model:vvo}. \eqref{eq:vvo_ref}
fixes the slack bus voltage angle. \eqref{eq:vvo_pg}--\eqref{eq:vvo_vm} enforce operating limits, and
\eqref{eq:vvo_cb}--\eqref{eq:vvo_tap_budget} impose device integrality and
switching budgets. With \(\lamc=1\), \(\lamv=\lamq=\lamp=0\), and fixed
\(\CDEV=0\), \(\TDEV=\TDEV^0\), Model~\ref{model:vvo} reduces to the canonical ACOPF problem, 
which is considered as one input formulation.

\subsection{DC Optimal Power Flow}

DCOPF minimizes generation cost under linearized AC constraint approximations: voltage
magnitudes are fixed at \(\vm_i=1\) p.u., reactive power is ignored, and
\eqref{eq:dcopf_balance} enforces linearized active power balance. Equations~\eqref{eq:dcopf_flow},
\eqref{eq:dcopf_thermal}, \eqref{eq:dcopf_pg}, and~\eqref{eq:dcopf_ang}
impose DC branch flow, thermal, generator, and angle-difference restrictions.

\modeltitle{model:dcopf}{DC Optimal Power Flow}
\begin{align}
  \min_{\PG,\VA}\quad &\sum_{i\in\GENERATORS}\cost_{i}(\PG_{i}) \label{eq:dcopf_obj}\\
  \text{s.t.}\quad
  &\thetaref=0, \nonumber\\
  &\sum_{e\in\Eiplus}\PF_{e}-\sum_{e\in\Eiminus}\PT_{e}
    =\PG_{i}-\pd_{i}, \quad\forall\,i\in\NODES, \label{eq:dcopf_balance}\\
  &\PF_{ij}=-\bij(\VA_{i}-\VA_{j}), \quad\forall\,ij\in\EDGES, \label{eq:dcopf_flow}\\
  &|\PF_{ij}|\leq\smax_{ij}, \quad\forall\,ij\in\EDGES, \label{eq:dcopf_thermal}\\
  &\pgmin_{i}\leq\PG_{i}\leq\pgmax_{i}, \quad\forall\,i\in\GENERATORS, \label{eq:dcopf_pg}\\
  &\dvamin_{ij}\leq\VA_{i}-\VA_{j}\leq\dvamax_{ij},
    \quad\forall\,ij\in\EDGES. \label{eq:dcopf_ang}
\end{align}
\modelbottomrule


\subsection{AC Power Flow}

Given network parameters, loads, bus-type specifications, active generation
setpoints, and fixed device positions, an AC power flow (ACPF) fixes the
non-slack active generations and the reference angle, then solves~\eqref{eq:kcl}--\eqref{eq:ohm:complex} to recover the
corresponding voltage, reactive generation, and branch flows. This paper follows the
standard PowerModels.jl implementation~\cite{PowerModels}. When ACPF converges,
the solution satisfies~\eqref{eq:kcl}--\eqref{eq:ohm:complex} up to solver tolerance, but it may still violate operating
limits such as~\eqref{eq:thermal}, \eqref{eq:vvo_qg}, and \eqref{eq:vvo_vm}.

\section{Methodology}
\label{sec:methodology}

\subsection{Homotopy-Guided Rounding}
\label{sec:homotopy_guided_rounding}

Model~\ref{model:vvo} is a nonconvex MINLP. Although it can be passed directly
to a mixed-integer solver such as Gurobi~\cite{gurobi}, prior work found off-the-shelf
solvers intractable in practice and introduced a relax--round--resolve
heuristic for transmission VVO~\cite{tong2026vvo}. The baseline heuristic first solves a
continuous relaxation, rounds fractional device positions, and then re-solves with
the integer device positions fixed. A key numerical difficulty is that 
rounding can create a \emph{rounding shock}: a solution feasible for fractional device values 
may not be feasible after those values are rounded. This paper replaces the one-shot re-solve with 
a homotopy continuation stage that moves gradually from relaxed to final integer device settings, 
avoiding the need for the continuous variables to absorb the full rounding effect in one step.

Define the continuous AC power-flow solution as
\begin{equation}
  \bv{x}=(\VM,\VA,\SG,\SF,\ST).
  \label{eq:state_vector}
\end{equation}
Let \(\DEV=(\CDEV,\TDEV)\) denote the variables controlling device positions.
The admissible CB and OLTC position sets are \(\CBPOSSET\) and
\(\TAPPOSSET\), and their maps into the AC power flow equations are
given by~\eqref{eq:cb-map} and~\eqref{eq:tap-map}.

Let
\(\DEVSET={\CBPOSSET}^{|\CBLOC|}\times{\TAPPOSSET}^{|\TAPLOC|}\)
define the set of admissible device positions, where \(\times\) denotes the
Cartesian product. The operator \(\operatorname{conv}(\cdot)\) denotes convex
hull. The relaxation stage relaxes integer constraint \(\DEV\in\DEVSET\) by
\begin{equation}
  \DEV\in\operatorname{conv}(\DEVSET)
  =
  {\operatorname{conv}(\CBPOSSET)}^{|\CBLOC|}
  \times
  {\operatorname{conv}(\TAPPOSSET)}^{|\TAPLOC|}.
  \label{eq:relaxed_device_set}
\end{equation}

Let \((\bv{x}^R,\DEV^R)\) denote the relaxed solution, where
\(\DEV^R=(\CDEV^R,\TDEV^R)\). The integer device setpoint is obtained by
\begin{equation}
  \DEV^Z
  =
  \Pi_{\DEVSET}(\DEV^R)
  =
  \left(
    \Pi_{\CBPOSSET}(\CDEV^R),
    \Pi_{\TAPPOSSET}(\TDEV^R)
  \right).
  \label{eq:device_projection}
\end{equation}
Here \(\Pi\) denotes nearest-neighbor projection onto the admissible integer settings.

Instead of re-solving the fixed stage immediately at
\(\DEV=\DEV^Z\), the homotopy continuation stage introduces a scalar continuation
parameter \(\alpha\in[0,1]\) and moves the device positions from their relaxed
values to their rounded values along the path
\begin{equation}
  \DEV(\alpha)=(1-\alpha)\DEV^R+\alpha\DEV^Z,
  \qquad \alpha\in[0,1].
  \label{eq:homotopy_path}
\end{equation}
At each \(\alpha\), \(\DEV(\alpha)\) is fixed, which also fixes the induced
susceptance and tap ratios, while the solver reoptimizes the continuous
network variables. Intermediate values of \(\CDEV_i(\alpha)\) and
\(\TDEV_{ij}(\alpha)\) are numerical continuation points that are not necessarily implementable settings.

For each \(\alpha\), \(\CB(\alpha)\) and \(\TAP(\alpha)\) affect the network equations
through~\eqref{eq:shunt} and~\eqref{eq:Ybr}. Let \(F\) collect the equality
constraints~\eqref{eq:kcl}, \eqref{eq:ohm:complex}, and~\eqref{eq:vvo_ref}, and
let \(h\) collect the inequality constraints~\eqref{eq:thermal}
and~\eqref{eq:vvo_pg}--\eqref{eq:vvo_vm}. The fixed-parameter homotopy subproblem for each \(\alpha\) is
\modeltitle{model:homotopy}{Fixed-Parameter Homotopy Subproblem}
\begin{align}
  \min_{\bv{x}}\quad
  & \psi(\VM,\QG,\PG) \label{eq:homotopy_obj}\\
  \text{s.t.}\quad
  & F\!\left(\bv{x};
      \bs+\CB(\alpha),
      \TAP(\alpha)\right)=0, \label{eq:homotopy_eq}\\
  & h\!\left(\bv{x};
      \bs+\CB(\alpha),
      \TAP(\alpha)\right)\le 0. \label{eq:homotopy_ineq}
\end{align}
\modelhomotopybottomrule

The continuation is solved over a sequence
\[
  0=\alpha_0 < \alpha_1 < \cdots < \alpha_m=1.
\]
Here \(m\) is the number of homotopy steps, and each subproblem is warm-started from
the previous subproblem's primal and dual solution. The final solution at \(\alpha=1\) 
is then checked for feasibility against the constraints of Model~\ref{model:vvo}.

\subsection{Input Formulation}

Model~\ref{model:vvo} is evaluated under the two workflows in
Figure~\ref{fig:workflow}. Let \(\XACOPF\) and \(\XACPF\) denote the ACOPF and
ACPF solutions, each with the components in~\eqref{eq:state_vector}. Let \(\PG(\bv{x})\)
denote the active power generation component of an AC solution vector \(\bv{x}\), and
let \(\XVVO=(\bv{x},\CDEV,\TDEV)\) collect the VVO variables. In
ACOPF$\to$VVO, VVO is warm-started from \(\XACOPF\), with
\(\pgref=\PG(\XACOPF)\), and the resulting solution is denoted
\(\XVVOACOPF\). In DCOPF$\to$ACPF$\to$VVO, the DCOPF solution
\((\PG,\PF,\VA)\) is used to prescribe the ACPF input. Since DCOPF does
not provide voltage magnitudes, this paper assumes \(\VM_i=1\) p.u. 
and \(\PT_{ij}=-\PF_{ij}\) for each bus and branch respectively to initialize the ACPF solver. The ACPF solution
then warm-starts VVO with \(\pgref=\PG(\XACPF)\), and the resulting VVO solution is denoted
\(\XVVOACPF\). 

\usetikzlibrary{fit,positioning,shapes.geometric,arrows.meta,backgrounds,matrix}


\begin{figure}[t]
\centering
\begin{tikzpicture}[
  proc/.style={rectangle, rounded corners=3pt, draw, thick,
               minimum height=0.58cm, minimum width=1.08cm,
               align=center, font=\footnotesize},
  greenproc/.style={proc, fill=green!12, draw=green!60!black},
  blueproc/.style={proc, fill=blue!10, draw=blue!50!black},
  grayproc/.style={proc, fill=gray!12, draw=gray!60},
  varbox/.style={rectangle, rounded corners=2pt, draw=gray!50,
                 fill=gray!6, inner sep=2.5pt, align=center, font=\scriptsize},
  arr/.style={->, >=Stealth, thick, gray!70},
  downarr/.style={->, >=Stealth, thick, gray!50},
  node distance=0.32cm and 0.48cm,
]

\node[grayproc] (caseAcpf) {Case\\file};
\node[greenproc, right=of caseAcpf] (dcopf) {DCOPF};
\node[greenproc, right=of dcopf] (acpf) {ACPF};
\node[blueproc, right=of acpf] (vvoAcpf) {VVO};

\draw[arr] (caseAcpf) -- (dcopf);
\draw[arr] (dcopf) -- (acpf);
\draw[arr] (acpf) -- (vvoAcpf);

\node[varbox, below=0.45cm of dcopf] (dcopfBox) {%
  $\PG$\\[-1pt]$\PF$\\[-1pt]$\VA$%
};
\node[varbox, below=0.45cm of acpf] (acpfBox) {$\XACPF$};
\node[varbox, below=0.45cm of vvoAcpf] (vvoAcpfBox) {$\XVVOACPF$};

\draw[downarr] (dcopf.south) -- (dcopfBox.north);
\draw[downarr] (acpf.south) -- (acpfBox.north);
\draw[downarr] (vvoAcpf.south) -- (vvoAcpfBox.north);

\node[grayproc, below=1.78cm of caseAcpf] (caseAcopf) {Case\\file};
\node[greenproc] (acopf) at (acpf |- caseAcopf) {ACOPF};
\node[blueproc] (vvoAcopf) at (vvoAcpf |- caseAcopf) {VVO};

\draw[arr] (caseAcopf) -- (acopf);
\draw[arr] (acopf) -- (vvoAcopf);

\node[varbox, below=0.45cm of acopf] (acopfBox) {$\XACOPF$};
\node[varbox, below=0.45cm of vvoAcopf] (vvoAcopfBox) {$\XVVOACOPF$};

\draw[downarr] (acopf.south) -- (acopfBox.north);
\draw[downarr] (vvoAcopf.south) -- (vvoAcopfBox.north);

\end{tikzpicture}
\caption{ACOPF and DCOPF$\to$ACPF input workflows for VVO.}
\label{fig:workflow}
\end{figure}
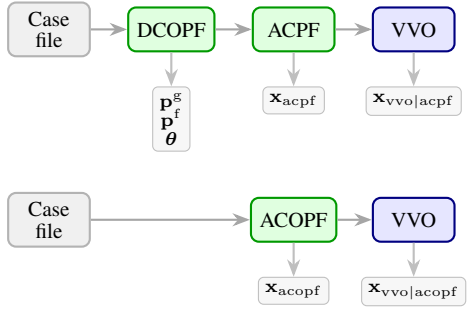

\subsection{Solution Method}
\label{sec:solution_method}

Algorithm~\ref{alg:homotopy_device_parameter_vvo} implements the construction
above using Ipopt~\cite{ipopt} for the relaxed VVO and homotopy continuation nonlinear program (NLP) solves.
If the projection in~\eqref{eq:device_projection} violates the budget in~\eqref{eq:vvo_cb_budget}--\eqref{eq:vvo_tap_budget},
the rounded point is repaired greedily by moving one device at a time toward its reference
setting while staying closest to the relaxed fractional value. The homotopy
stage then solves Model~\ref{model:homotopy} along~\eqref{eq:homotopy_path},
warm-starting each subproblem from the previous primal and dual solutions and
adapting the step size using the NLP solver's iteration counts.

\begin{algorithm}[!t]
\caption{Homotopy VVO Heuristic}
\label{alg:homotopy_device_parameter_vvo}
\begin{algorithmic}[1]
\State \textbf{Input:} Network data; setpoint
  \(\bv{x}^{0}\in\{\XACOPF,\XACPF\}\); device settings
  \(\CBPOSSET,\TAPPOSSET\); budgets \(B_{\beta},B_{\tau}\);
  step controls \(\Delta\alpha_0,\Delta\alpha_{\min}\)
\Statex
\State \emph{Relaxed VVO}
\State Relax \(\DEV\in\DEVSET\) to
  \(\DEV\in\operatorname{conv}(\DEVSET)\)
\State Solve Model~\ref{model:vvo} from \(\bv{x}^{0}\)
\State Retrieve primal and dual solutions
  \((\bv{x}^{R},\DEV^{R},\vect{\lambda}^{R})\)
\Statex
\State \emph{Round device settings}
\For{each capacitor bank bus \(i\in\CBLOC\)}
  \State \(\CDEV_i^{Z}\gets\Pi_{\CBPOSSET}(\CDEV_i^{R})\)
\EndFor
\For{each OLTC branch \(ij\in\TAPLOC\)}
  \State \(\TDEV_{ij}^{Z}\gets\Pi_{\TAPPOSSET}(\TDEV_{ij}^{R})\)
\EndFor
\State 
  while satisfying~\eqref{eq:vvo_cb_budget}--\eqref{eq:vvo_tap_budget}
\Statex
\State \emph{Homotopy continuation}
\State \(\alpha\gets0\), \(\Delta\alpha\gets\Delta\alpha_0\),
  \((\bv{x},\vect{\lambda})\gets(\bv{x}^{R},\vect{\lambda}^{R})\)
\While{\(\alpha<1\)}
  \State \(\alpha^{+}\gets\min\{1,\alpha+\Delta\alpha\}\)
  \State Fix \(\DEV\equiv(1-\alpha^{+})\DEV^{R}+\alpha^{+}\DEV^{Z}\)
  \State Solve Model~\ref{model:homotopy} warm-started from
    \((\bv{x},\vect{\lambda})\)
  \If{the subproblem is locally solved}
    \State Update \((\bv{x},\vect{\lambda})\) to the new solution
    \State \(\alpha\gets\alpha^{+}\)
    \State Adapt \(\Delta\alpha\) based on the NLP iteration count
  \Else
    \State \(\Delta\alpha\gets\Delta\alpha/2\)
    \If{\(\Delta\alpha<\Delta\alpha_{\min}\)}
      \State \Return no solution found
    \EndIf
  \EndIf
\EndWhile
\If{the final setpoint satisfies Model~\ref{model:vvo}}
  \State \Return \(\XVVO=(\bv{x},\DEV^Z)\)
\Else
  \State \Return no solution found
\EndIf
\end{algorithmic}
\end{algorithm}

Both input workflows use the same model and algorithm. Only the initial setpoint differs,
isolating the effect of input quality within each experimental setting.

\section{Numerical Experiments}\label{sec:experiments}

\subsection{Experimental Setup}

Experiments consider 13 PGLib cases from IEEE, PEGASE, and RTE systems
ranging from 118 to 13,659 buses \cite{pglib}.
Two device range settings are reported: transformer taps can move either
\(\pm3\) steps, or over the full \(\pm16\) range from reference, corresponding to
\(\TAPPOSSET=\{-3,-2,\ldots,2,3\}\) and
\(\TAPPOSSET=\{-16,-15,\ldots,15,16\}\), respectively.
Capacitor bank settings are the same: \(\CBPOSSET=\{-1,0,1,2\}\), which corresponds to 0, 1, 2, or 3 
active modules assuming one module is active at the reference susceptance value.
The device movement budget multipliers are set to
\(\rho_{\tau}=\rho_{\beta}=1\), allowing one
tap or CB module switch on average per OLTC branch and per CB bus.
To balance Volt/VAR control against redispatch and generation cost, the objective
weights are \(\lamp=\lamv=\lamq=\lamc=1\) for all experiments.
Unless otherwise stated, the default voltage and reactive power references are set
to \(\vref_i=1\) p.u. and \(\qgref_i=0\).
In practice, operators can replace these values with system-specific targets. Taking
\(\lamp \rightarrow \infty\), or equivalently constraining
\(\PG_i=\pgref_i\) for all \(i\neq\mathrm{slack}\), would forbid redispatch
other than at the slack bus.

All models are implemented in Julia~1.11 with JuMP~\cite{JuMP} and
PowerModels.jl~\cite{PowerModels}, solved with Ipopt on the PACE Phoenix
cluster~\cite{PACE} configured with Intel Xeon 6226 @ 2.7\,GHz, 64\,GB RAM, 8 CPU cores,
and a 4\,h time limit unless otherwise noted.

\subsection{Violations Across Test Cases}

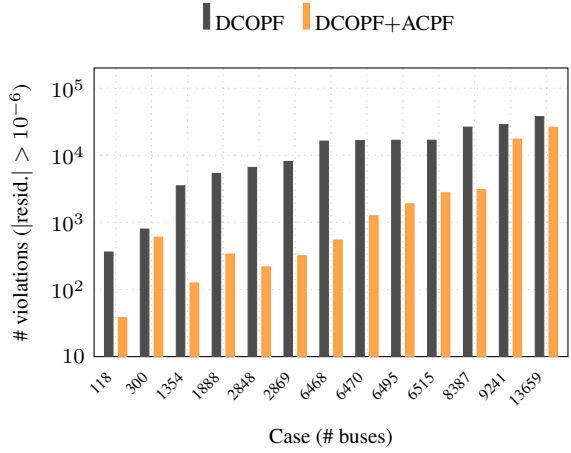
\begin{figure}[t]
\centering
\begin{tikzpicture}
\begin{axis}[
  ybar,
  ymode=log,
  log origin=infty,
  bar width=3.2pt,
  width=\columnwidth,
  height=5.4cm,
  enlarge x limits=0.05,
  ylabel={\footnotesize \# violations ($|\text{resid.}|>10^{-6}$)},
  ylabel style={font=\footnotesize},
  ymin=10, ymax=200000,
  ytick={10,100,1000,10000,100000},
  yticklabels={$10$,$10^2$,$10^3$,$10^4$,$10^5$},
  yticklabel style={font=\footnotesize},
  xtick=data,
  xticklabels={118,300,1354,1888,2848,2869,6468,6470,6495,6515,8387,9241,13659},
  xticklabel style={font=\tiny, rotate=45, anchor=east},
  xlabel={\footnotesize Case (\# buses)},
  xlabel style={font=\footnotesize},
  legend style={font=\footnotesize, at={(0.5,1.08)}, anchor=south,
                legend columns=2, draw=none,
                /tikz/every even column/.append style={column sep=0.18cm}},
  legend image code/.code={%
    \draw[#1] (0cm,-0.06cm) rectangle (0.08cm,0.28cm);
  },
  legend cell align=left,
  grid=major,
  grid style={dotted, gray!40},
  tick style={draw=none},
]
\addplot[fill=black!70, draw=black!70] coordinates {
  (1,  362)
  (2,  796)
  (3,  3508)
  (4,  5353)
  (5,  6604)
  (6,  8099)
  (7,  16280)
  (8,  16588)
  (9,  16737)
  (10, 16756)
  (11, 26269)
  (12, 28657)
  (13, 37681)
};
\addlegendentry{DCOPF}

\addplot[fill=orange!70, draw=orange!80] coordinates {
  (1,  38)
  (2,  604)
  (3,  125)
  (4,  337)
  (5,  217)
  (6,  318)
  (7,  546)
  (8,  1255)
  (9,  1888)
  (10, 2750)
  (11, 3094)
  (12, 17425)
  (13, 26141)
};
\addlegendentry{DCOPF$+$ACPF}

\end{axis}
\end{tikzpicture}
\caption{DCOPF and DCOPF$+$ACPF constraint violations with $|\text{residual}|>10^{-6}$.}
\label{fig:violations}
\end{figure}

Figure~\ref{fig:violations} reports the number of AC constraint violations larger than 
$10^{-6}$ for the DCOPF and DCOPF$+$ACPF setpoints. Since the DCOPF solution lacks AC variables, \(\VM_i=1, \QG_i=0 \;\forall i\), \(\QF_{ij}\!=\!\QT_{ij}=0, \PT_{ij}=-\PF_{ij} \;\forall ij\) are assumed.
Depending on system size, DCOPF solutions produce 362--37,681 violations,
primarily in nodal power balance~\eqref{eq:kcl} and branch
flow~\eqref{eq:ohm:complex} equations.

Applying ACPF substantially reduces AC equation violations in count and
magnitude, but does not enforce all operating limits in Model~\ref{model:vvo}.
Reactive generation and thermal limit violations remain, with 38--26,141 total
violations across the cases and residuals up to \(10^{-2}\) p.u. These
remaining violations show that ACPF is not a complete feasibility recovery
layer. Coordinated control is required, motivating a downstream VVO framework.


\subsection{Solution Quality}

Tables~\ref{tab:quality} and~\ref{tab:quality-full} compare four solution
classes. The base ACOPF row reports \(\XACOPF\), obtained by solving the canonical
ACOPF problem with CB offsets fixed at \(\CDEV_i=0\) and tap positions fixed at
\(\TDEV_{ij}=\TDEV^0_{ij}\). The proposed homotopy-guided heuristic solves Model~\ref{model:vvo} using
Algorithm~\ref{alg:homotopy_device_parameter_vvo} under the two input workflows,
yielding \(\XVVOACOPF\) and \(\XVVOACPF\). The vanilla VVO heuristic from~\cite{tong2026vvo}
is initialized from \(\XACOPF\) and uses the same objective coefficients as the proposed heuristic. It
also optimizes \(\CDEV_i\) and \(\TDEV_{ij}\), but omits the movement budgets~\eqref{eq:vvo_cb_budget}--\eqref{eq:vvo_tap_budget}
and uses a one-shot re-solve after rounding.
Infeasible runs are marked with dashes.
To keep the tables compact, eight representative cases are reported. Cases 1354, 2848, 2869, 6470, and 6495 are omitted
because their system sizes are already represented by the selected cases.

The reported quality metrics are
\begin{align}
  \mathrm{MAE}_v
  &= \frac{1}{|\NODES|}\sum_{i\in\NODES}|\VM_i-\vref_i|,
  \label{eq:vmae}\\
  \mathrm{MAE}_q
  &= \frac{1}{|\GENERATORS|}\sum_{i\in\GENERATORS}|\QG_i-\qgref_i|,
  \label{eq:qmae}\\
  \mathrm{MAE}_p
  &= \frac{1}{|\GENERATORS|}\sum_{i\in\GENERATORS}|\PG_i-\pgref_i|,
  \label{eq:pmae}\\
  \%\Delta c
  &=100\,\frac{\psi_c(\PG)-\psi_c(\PG(\XACOPF))}
    {\psi_c(\PG(\XACOPF))} .
  \label{eq:costdelta}
\end{align}
where \(\PG(\XACOPF)\) is the ACOPF dispatch. Lower
\(\mathrm{MAE}_v\), \(\mathrm{MAE}_q\), and \(\mathrm{MAE}_p\) indicate smaller
voltage deviations, improved reactive-power regulation, and smaller redispatch from \(\pgref\),
respectively; negative \(\%\Delta c\) indicates cost reduction.
All runtimes are reported in seconds.
\(T_r\) is the relaxed VVO solve time. For Vanilla VVO, \(T_f\) is the one-shot 
re-solve time after rounding. For homotopy-guided VVO, \(T_f\) is the total
time in the homotopy continuation stage, including all intermediate subproblems and the terminal solve at \(\alpha=1\).
For base ACOPF rows, \(T_r\) and \(T_f\) are not applicable.

\textbf{Feasibility.}
Across all 52 runs, covering 13 cases, two device configurations, and two input workflows,
the homotopy-guided heuristic succeeds in every instance. In contrast, vanilla VVO succeeds
in 42/52 runs. For both ACOPF and DCOPF$\to$ACPF inputs, vanilla VVO fails on IEEE 300 with
tap $\pm3$, and on IEEE 300, RTE 1888, 6495, and 6515 with tap $\pm16$.

\textbf{Input sensitivity.}
Homotopy-guided VVO does not require the input setpoint to be feasible for
Model~\ref{model:vvo}. Starting from \(\XACPF\), which may violate operating
limits, \(\XVVOACPF\) restores Model~\ref{model:vvo} feasibility while
delivering solution quality close to \(\XVVOACOPF\). The main difference appears in \(\mathrm{MAE}_p\): \(\XVVOACOPF\) uses \(\pgref=\PG(\XACOPF)\), whereas \(\XVVOACPF\) uses \(\pgref=\PG(\XACPF)\). The larger \(\mathrm{MAE}_p\) for \(\XVVOACPF\) indicates that active generation must move
farther from the input dispatch to restore feasibility. In contrast, other quality metrics
remain similar between \(\XVVOACOPF\) and \(\XVVOACPF\) across device ranges. Relative to the base ACOPF setpoints, VVO
setpoints can increase voltage deviations, but the resulting voltages remain within secure limits, while reactive dispatch and
generation cost decrease in nearly every reported row.


\begin{table}[t]
\centering
\scriptsize
\setlength{\tabcolsep}{2pt}
\caption{Solution quality for tap $\pm3$, CB 0--3.}
\label{tab:quality}
\begin{tabular}{l l rrr r rr}
\toprule
\textbf{Case} & \textbf{Input}
  & $\mathrm{MAE}_{v}$
  & $\mathrm{MAE}_{q}$
  & $\mathrm{MAE}_{p}$
  & $\%\!\Delta c$
  & $T_r$
  & $T_f$ \\
\midrule
\multirow{4}{*}{118}
  & base              & 0.034 & 38.00 & 0.00 & 0.00 & $-$ & $-$ \\
  & Vanilla VVO       & \textbf{0.036} & \textbf{35.51} & 0.25 & $\mathbf{-0.06}$ & \textbf{0.2} & \textbf{0.2} \\
  & ACOPF             & \textbf{0.036} & 37.05 & \textbf{0.19} & $-$0.05 & 3.8 & 5.1 \\
  & DCOPF$\!\to\!$ACPF & \textbf{0.036} & 37.10 & 5.76 & $-$0.05 & 3.6 & 5.1 \\
\midrule
\multirow{4}{*}{300}
  & base              & 0.030 & 126.78 & 0.00 & 0.00 & $-$ & $-$ \\
  & Vanilla VVO       & $-$ & $-$ & $-$ & $-$ & $-$ & $-$ \\
  & ACOPF             & \textbf{0.031} & \textbf{129.58} & \textbf{29.48} & $\mathbf{-2.29}$ & 4.6 & \textbf{6.5} \\
  & DCOPF$\!\to\!$ACPF & \textbf{0.031} & 130.05 & 44.35 & $\mathbf{-2.29}$ & \textbf{4.4} & \textbf{6.5} \\
\midrule
\multirow{4}{*}{1888}
  & base              & 0.040 & 30.02 & 0.00 & 0.00 & $-$ & $-$ \\
  & Vanilla VVO       & \textbf{0.070} & \textbf{24.52} & 28.46 & $\mathbf{-2.19}$ & \textbf{8.9} & \textbf{4.9} \\
  & ACOPF             & 0.072 & 24.91 & 28.47 & $\mathbf{-2.19}$ & 15.3 & 32.5 \\
  & DCOPF$\!\to\!$ACPF & 0.072 & 24.92 & \textbf{5.05} & $\mathbf{-2.19}$ & 13.2 & 32.3 \\
\midrule
\multirow{4}{*}{6468}
  & base              & 0.048 & 34.30 & 0.00 & 0.00 & $-$ & $-$ \\
  & Vanilla VVO       & \textbf{0.059} & 28.25 & 0.47 & $-$0.13 & \textbf{49.9} & \textbf{31.0} \\
  & ACOPF             & 0.061 & \textbf{27.87} & \textbf{0.32} & $\mathbf{-0.14}$ & 50.1 & 162.3 \\
  & DCOPF$\!\to\!$ACPF & 0.061 & 27.98 & 10.91 & $\mathbf{-0.14}$ & 61.2 & 140.1 \\
\midrule
\multirow{4}{*}{6515}
  & base              & 0.047 & 40.03 & 0.00 & 0.00 & $-$ & $-$ \\
  & Vanilla VVO       & \textbf{0.054} & \textbf{32.03} & \textbf{0.54} & $-$0.14 & \textbf{30.8} & \textbf{23.3} \\
  & ACOPF             & 0.059 & 32.63 & 0.80 & $\mathbf{-0.20}$ & 48.2 & 138.0 \\
  & DCOPF$\!\to\!$ACPF & 0.059 & 32.63 & 9.43 & $\mathbf{-0.20}$ & 120.9 & 126.1 \\
\midrule
\multirow{4}{*}{8387}
  & base              & 0.079 & 95.52 & 0.00 & 0.00 & $-$ & $-$ \\
  & Vanilla VVO       & \textbf{0.080} & \textbf{79.43} & 6.40 & $-$1.09 & \textbf{64.1} & \textbf{107.5} \\
  & ACOPF             & 0.083 & 80.37 & \textbf{5.80} & $\mathbf{-1.11}$ & 96.2 & 339.2 \\
  & DCOPF$\!\to\!$ACPF & 0.083 & 80.26 & 38.92 & $-$1.09 & 412.6 & 342.6 \\
\midrule
\multirow{4}{*}{9241}
  & base              & 0.060 & 62.99 & 0.00 & 0.00 & $-$ & $-$ \\
  & Vanilla VVO       & \textbf{0.063} & 32.02 & 1.47 & $-$0.07 & 522.4 & \textbf{112.5} \\
  & ACOPF             & 0.081 & \textbf{26.95} & \textbf{1.35} & $\mathbf{-0.24}$ & 707.9 & 299.8 \\
  & DCOPF$\!\to\!$ACPF & 0.081 & 27.15 & 10.95 & $\mathbf{-0.24}$ & \textbf{510.2} & 296.4 \\
\midrule
\multirow{4}{*}{13659}
  & base              & 0.068 & 28.83 & 0.00 & 0.00 & $-$ & $-$ \\
  & Vanilla VVO       & \textbf{0.067} & 15.26 & 0.66 & $-$0.05 & 762.4 & \textbf{87.9} \\
  & ACOPF             & 0.084 & \textbf{13.56} & \textbf{0.45} & $\mathbf{-0.14}$ & \textbf{509.4} & 440.8 \\
  & DCOPF$\!\to\!$ACPF & 0.084 & 13.58 & 3.26 & $\mathbf{-0.14}$ & 565.6 & 448.4 \\
\bottomrule
\end{tabular}
\end{table}


\begin{table}[t]
\centering
\scriptsize
\setlength{\tabcolsep}{2pt}
\caption{Solution quality for tap $\pm16$, CB 0--3.}
\label{tab:quality-full}
\begin{tabular}{l l rrr r rr}
\toprule
\textbf{Case} & \textbf{Input}
  & $\mathrm{MAE}_{v}$
  & $\mathrm{MAE}_{q}$
  & $\mathrm{MAE}_{p}$
  & $\%\!\Delta c$
  & $T_r$
  & $T_f$ \\
\midrule
\multirow{4}{*}{118}
  & base              & 0.034 & 38.00 & 0.00 & 0.00 & $-$ & $-$ \\
  & Vanilla VVO       & 0.037 & \textbf{34.09} & 0.41 & $\mathbf{-0.07}$ & \textbf{0.2} & \textbf{0.1} \\
  & ACOPF             & \textbf{0.036} & 36.24 & \textbf{0.08} & $-$0.04 & 3.9 & 5.1 \\
  & DCOPF$\!\to\!$ACPF & \textbf{0.036} & 36.29 & 5.88 & $-$0.04 & 3.8 & 5.3 \\
\midrule
\multirow{4}{*}{300}
  & base              & 0.030 & 126.78 & 0.00 & 0.00 & $-$ & $-$ \\
  & Vanilla VVO       & $-$ & $-$ & $-$ & $-$ & $-$ & $-$ \\
  & ACOPF             & \textbf{0.031} & \textbf{130.97} & \textbf{35.56} & $\mathbf{-2.87}$ & \textbf{4.5} & \textbf{6.4} \\
  & DCOPF$\!\to\!$ACPF & \textbf{0.031} & 131.45 & 39.35 & $\mathbf{-2.87}$ & 4.7 & 7.0 \\
\midrule
\multirow{4}{*}{1888}
  & base              & 0.040 & 30.02 & 0.00 & 0.00 & $-$ & $-$ \\
  & Vanilla VVO       & $-$ & $-$ & $-$ & $-$ & $-$ & $-$ \\
  & ACOPF             & \textbf{0.072} & \textbf{24.00} & 28.48 & $\mathbf{-2.21}$ & 19.4 & 32.8 \\
  & DCOPF$\!\to\!$ACPF & \textbf{0.072} & 24.01 & \textbf{5.02} & $\mathbf{-2.21}$ & \textbf{16.6} & \textbf{32.0} \\
\midrule
\multirow{4}{*}{6468}
  & base              & 0.048 & 34.30 & 0.00 & 0.00 & $-$ & $-$ \\
  & Vanilla VVO       & \textbf{0.058} & 32.08 & 6.81 & 0.27 & 604.0 & \textbf{36.2} \\
  & ACOPF             & 0.062 & \textbf{25.97} & \textbf{1.09} & $\mathbf{-0.19}$ & \textbf{51.2} & 155.5 \\
  & DCOPF$\!\to\!$ACPF & 0.062 & 25.99 & 10.90 & $\mathbf{-0.19}$ & 59.4 & 154.0 \\
\midrule
\multirow{4}{*}{6515}
  & base              & 0.047 & 40.03 & 0.00 & 0.00 & $-$ & $-$ \\
  & Vanilla VVO       & $-$ & $-$ & $-$ & $-$ & $-$ & $-$ \\
  & ACOPF             & \textbf{0.060} & \textbf{29.63} & \textbf{2.22} & $\mathbf{-0.32}$ & \textbf{44.0} & \textbf{132.4} \\
  & DCOPF$\!\to\!$ACPF & \textbf{0.060} & 29.75 & 9.57 & $\mathbf{-0.32}$ & 89.7 & 138.4 \\
\midrule
\multirow{4}{*}{8387}
  & base              & 0.079 & 95.52 & 0.00 & 0.00 & $-$ & $-$ \\
  & Vanilla VVO       & \textbf{0.073} & 83.20 & 16.71 & $\mathbf{-3.65}$ & \textbf{113.7} & \textbf{87.0} \\
  & ACOPF             & 0.078 & \textbf{81.54} & \textbf{12.31} & $-$2.73 & 123.1 & 336.5 \\
  & DCOPF$\!\to\!$ACPF & 0.078 & 81.89 & 39.56 & $-$2.71 & 509.6 & 325.3 \\
\midrule
\multirow{4}{*}{9241}
  & base              & 0.060 & 62.99 & 0.00 & 0.00 & $-$ & $-$ \\
  & Vanilla VVO       & \textbf{0.065} & 33.05 & 2.05 & $-$0.10 & 2441.7 & \textbf{152.6} \\
  & ACOPF             & 0.083 & \textbf{26.01} & \textbf{1.60} & $\mathbf{-0.27}$ & \textbf{532.1} & 312.4 \\
  & DCOPF$\!\to\!$ACPF & 0.083 & 26.15 & 10.84 & $\mathbf{-0.27}$ & 557.8 & 298.0 \\
\midrule
\multirow{4}{*}{13659}
  & base              & 0.068 & 28.83 & 0.00 & 0.00 & $-$ & $-$ \\
  & Vanilla VVO       & \textbf{0.063} & 14.18 & 0.72 & $-$0.08 & 2397.1 & \textbf{100.3} \\
  & ACOPF             & 0.087 & 12.08 & \textbf{0.62} & $\mathbf{-0.18}$ & 700.5 & 429.3 \\
  & DCOPF$\!\to\!$ACPF & 0.087 & \textbf{12.02} & 3.34 & $\mathbf{-0.18}$ & \textbf{609.8} & 451.4 \\
\bottomrule
\end{tabular}
\end{table}

\textbf{Runtimes.}
Although the homotopy continuation stage adds intermediate subproblems, it continues to produce high
quality solutions and improves runtime stability relative to vanilla VVO, especially in the full
tap range setting where vanilla VVO can be slower on large cases.
Moreover, \(T_r\) and \(T_f\) remain similar for \(\XACOPF\) and
\(\XACPF\) inputs in most cases, indicating that restoring feasibility
from DCOPF$\to$ACPF inputs does not introduce a systematic runtime premium.

\textbf{Device movement.}
Figure~\ref{fig:device-movement} illustrates the practical effect of
the device movement budgets. Only cases where both approaches produce
feasible solutions are shown. Relative to unbudgeted vanilla VVO, the budgeted
homotopy runs substantially reduce tap movement, especially in the full-range
setting, while keeping CB movement close to or below the vanilla values and
maintaining similar solution quality.


\pgfplotsset{
  deviceMovementAxis/.style={
    /pgfplots/ybar,
    /pgfplots/bar width=2.1pt,
    /pgfplots/width=\linewidth,
    /pgfplots/height=4.35cm,
    /pgfplots/ymin=8,
    /pgfplots/enlarge x limits=0.03,
    /pgfplots/ylabel={\scriptsize total abs. steps},
    /pgfplots/ylabel style={font=\scriptsize},
    /pgfplots/xlabel={\scriptsize Case (\# buses)},
    /pgfplots/xlabel style={font=\scriptsize},
    /pgfplots/tick style={draw=none},
    /pgfplots/xticklabel style={font=\tiny, rotate=45, anchor=east},
    /pgfplots/ytick={10,100,1000,10000},
    /pgfplots/yticklabels={$10$,$10^2$,$10^3$,$10^4$},
    /pgfplots/yticklabel style={font=\scriptsize},
    /pgfplots/grid=major,
    /pgfplots/grid style={dotted, gray!40},
    /pgfplots/title style={font=\footnotesize},
  }
}

\begin{figure*}[t]
\centering
\begin{tikzpicture}[baseline]
  \draw[fill=orange!75, draw=orange!85!black] (0,0) rectangle (0.22,0.14);
  \node[anchor=west, font=\small] at (0.28,0.07) {Budgeted homotopy};
  \draw[fill=black!60, draw=black!70] (3.55,0) rectangle (3.77,0.14);
  \node[anchor=west, font=\small] at (3.83,0.07) {Vanilla VVO};
\end{tikzpicture}

\vspace{0.15em}

\begin{minipage}[t]{0.49\textwidth}
\centering
\begin{tikzpicture}
\begin{axis}[
  ymode=log,
  log origin=infty,
  deviceMovementAxis,
  title={tap movement, tap $\pm3$},
  ymax=50000,
  xtick={1,2,3,4,5,6,7,8,9,10,11,12},
  xticklabels={118,1354,1888,2848,2869,6468,6470,6495,6515,8387,9241,13659},
]
\addplot[fill=orange!75, draw=orange!85!black] coordinates {
  (1,11) (2,239) (3,553) (4,783) (5,531) (6,1572)
  (7,1579) (8,1597) (9,1615) (10,2043) (11,2252) (12,6641)
};
\addplot[fill=black!60, draw=black!70] coordinates {
  (1,25) (2,601) (3,1221) (4,1766) (5,1277) (6,3382)
  (7,3429) (8,3550) (9,3659) (10,5294) (11,4945) (12,16030)
};
\end{axis}
\end{tikzpicture}
\end{minipage}\hfill
\begin{minipage}[t]{0.49\textwidth}
\centering
\begin{tikzpicture}
\begin{axis}[
  ymode=log,
  log origin=infty,
  deviceMovementAxis,
  title={CB movement, tap $\pm3$},
  ylabel={},
  ymax=15000,
  xtick={1,2,3,4,5,6,7,8,9,10,11,12},
  xticklabels={118,1354,1888,2848,2869,6468,6470,6495,6515,8387,9241,13659},
]
\addplot[fill=orange!75, draw=orange!85!black] coordinates {
  (1,14) (2,1065) (3,45) (4,46) (5,2167) (6,96)
  (7,73) (8,99) (9,101) (10,485) (11,7085) (12,8549)
};
\addplot[fill=black!60, draw=black!70] coordinates {
  (1,18) (2,1273) (3,53) (4,46) (5,2535) (6,130)
  (7,97) (8,121) (9,142) (10,716) (11,7052) (12,8925)
};
\end{axis}
\end{tikzpicture}
\end{minipage}

\vspace{0.35em}

\begin{minipage}[t]{0.49\textwidth}
\centering
\begin{tikzpicture}
\begin{axis}[
  ymode=log,
  log origin=infty,
  deviceMovementAxis,
  title={tap movement, tap $\pm16$},
  ymax=50000,
  xtick={1,2,3,4,5,6,7,8,9},
  xticklabels={118,1354,2848,2869,6468,6470,8387,9241,13659},
]
\addplot[fill=orange!75, draw=orange!85!black] coordinates {
  (1,10) (2,238) (3,779) (4,531) (5,1573)
  (6,1577) (7,2047) (8,2250) (9,6646)
};
\addplot[fill=black!60, draw=black!70] coordinates {
  (1,43) (2,1554) (3,3704) (4,3028) (5,6904)
  (6,7720) (7,18946) (8,11123) (9,34992)
};
\end{axis}
\end{tikzpicture}
\end{minipage}\hfill
\begin{minipage}[t]{0.49\textwidth}
\centering
\begin{tikzpicture}
\begin{axis}[
  ymode=log,
  log origin=infty,
  deviceMovementAxis,
  title={CB movement, tap $\pm16$},
  ylabel={},
  ymax=15000,
  xtick={1,2,3,4,5,6,7,8,9},
  xticklabels={118,1354,2848,2869,6468,6470,8387,9241,13659},
]
\addplot[fill=orange!75, draw=orange!85!black] coordinates {
  (1,14) (2,1079) (3,48) (4,2182) (5,96)
  (6,72) (7,487) (8,7086) (9,8548)
};
\addplot[fill=black!60, draw=black!70] coordinates {
  (1,18) (2,1379) (3,45) (4,2620) (5,130)
  (6,96) (7,720) (8,6944) (9,8997)
};
\end{axis}
\end{tikzpicture}
\end{minipage}
\caption{Total absolute device movements from reference positions for the homotopy and vanilla VVO methods.}
\label{fig:device-movement}
\end{figure*}
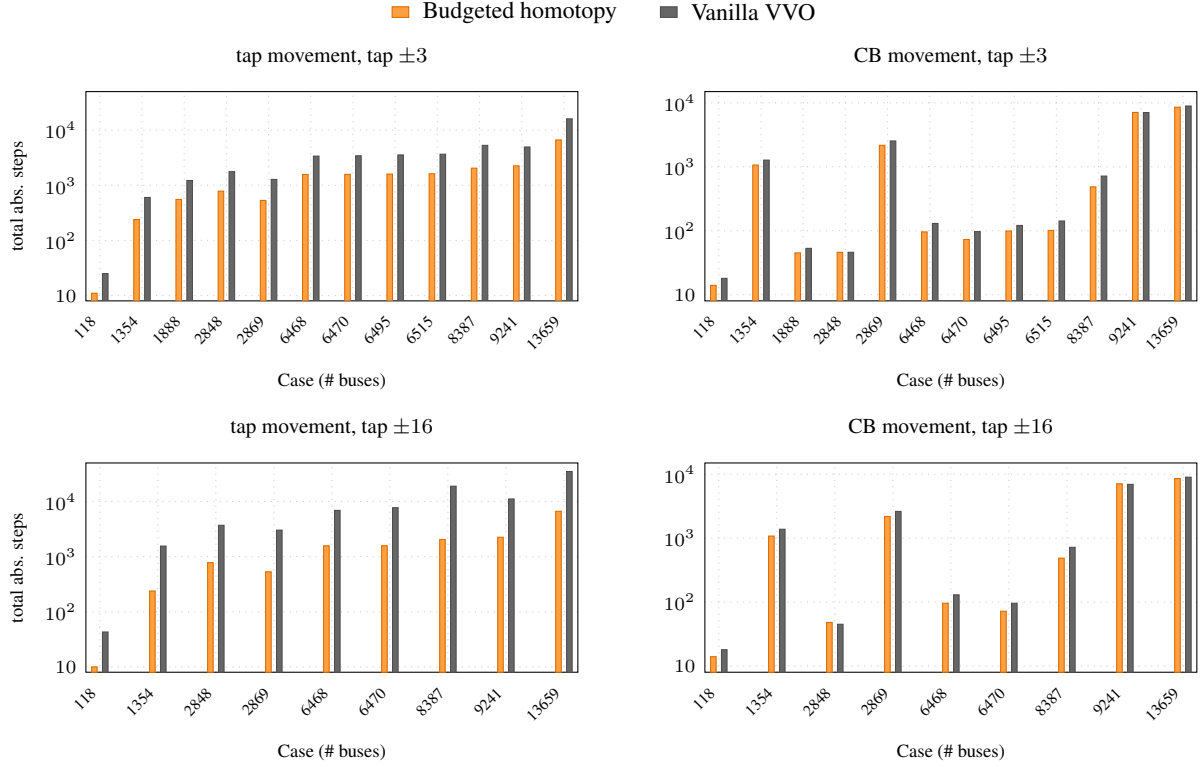

\subsection{Heuristic Baseline Comparison}

This paper compares homotopy-guided VVO with two baselines
under the DCOPF$\to$ACPF input workflow. The comparison is reported for tap
ranges $\pm3$ and $\pm16$ with 0--3 capacitor-bank modules, on the ten cases
where all three methods produce feasible solutions.

\textbf{Feasibility pump (FP-storm).}
The first baseline adapts the ``storm of feasibility pumps'' heuristic for
nonconvex MINLP~\cite{DAmbrosio2012StormFP}. FP-storm alternates between solving a mixed-integer (MILP) projection that enforces integer constraints and a relaxed NLP projection that enforces AC constraints. Given relaxed device reference (\(\CDEV^{\rm ref},\TDEV^{\rm ref}\)), the MILP projection selects admissible CB and transformer settings that are closest to the relaxed values, schematically minimizing
\begin{equation}
\sum_{i\in\CBLOC}|\CDEV_i-\CDEV_i^{\rm ref}|
+
\sum_{ij\in\TAPLOC}|\TDEV_{ij}-\TDEV_{ij}^{\rm ref}|,
\label{eq:fpstorm_milp}
\end{equation}
subject to the integer device sets, switching budgets, and anti-cycling constraints against revisiting previously tested configurations. The resulting integer solution (\(\CDEV^{\rm int},\TDEV^{\rm int}\)) is then tested by fixing the CB and OLTC
positions and solving the fixed-device VVO problem for the remaining
continuous AC variables as the NLP projection step. A feasible re-solve is accepted for quality evaluation. Otherwise, FP-storm solves a new NLP projection by minimizing the objective of~\eqref{eq:vvo_obj} plus a quadratic penalty toward the current integer target:
\begin{equation}
\psi(\VM,\QG,\PG)
+ \eta\lVert\CDEV-\CDEV^{\rm int}\rVert^2
+ \eta\lVert\TDEV-\TDEV^{\rm int}\rVert^2,
\label{eq:fpstorm_nlp_penalty}
\end{equation}
where \(\eta>0\) is a penalty weight. The NLP projection is subject to the AC
constraints and continuous relaxations of the device bounds. The fractional device values become
the next reference for the MILP projection.


\textbf{Sensitivity-guided rounding.}
The second baseline implements a merit function sensitivity approach \cite{Capitanescu2010SensitivityOPF}. 
At each iteration, one-step integer moves of \(\CDEV\) and \(\TDEV\) are ranked
by a merit score that estimates their first-order changes on the VVO objective~\eqref{eq:vvo_obj}
and operating limit residuals through the induced
shunt susceptance \(\CB_i(\CDEV_i)\) and tap ratio \(\TAP_{ij}(\TDEV_{ij})\)
in the AC network equations. After applying the selected move, the method
re-optimizes the continuous AC variables with fixed device settings. If this re-solve
fails, the method retains the same integer settings to solve a minimum-violation NLP that enforces AC equations
while minimizing operating limit slacks. The
procedure terminates when no improving move remains or a step limit is reached.

The sparse implementation evaluates these scores through linearized response
calculations. Let \(g(y,z)=0\) denote the equality system of
AC constraints with fixed device settings, where \(y\) collects the continuous AC
variables and \(z\) collects the discrete device controls. The method forms the
equality Jacobians \(J_y=\partial g/\partial y\) and
\(J_z=\partial g/\partial z\). For each candidate device direction \(j\), let
\(J_{z,j}\) denote the \(j\)th column of \(J_z\). The linearized response is
computed from
\begin{equation}
  J_y \Delta y_j = -J_{z,j}.
  \label{eq:sens-linear-response}
\end{equation}
Solving~\eqref{eq:sens-linear-response} gives \(\Delta y_j\), the first-order
change in \(y\) induced by a unit move in \(z_j\).
Let \(f(y)=\psi(\VM,\QG,\PG)\) denote the VVO objective
in~\eqref{eq:vvo_obj}, and let \(h(y,z)\leq 0\) collect the inequality
operating limits, including \eqref{eq:thermal} and~\eqref{eq:vvo_pg}--\eqref{eq:vvo_vm}.
For a signed one-step move \(\delta z_j\), the estimated objective and
inequality-residual changes are
\begin{equation}
  \begin{aligned}
    \widehat{\Delta f}_j
    &= \delta z_j\,\nabla_y f^\top \Delta y_j,\\
    \widehat{\Delta h}_j
    &= \delta z_j\,(K_y \Delta y_j + K_{z,j}),
  \end{aligned}
  \label{eq:sens-score-linearization}
\end{equation}
where \(\nabla_y f=\partial f/\partial y\), \(K_y=\partial h/\partial y\),
\(K_z=\partial h/\partial z\), and
\(K_{z,j}\) is the \(j\)th column of \(K_z\).
Equations~\eqref{eq:sens-linear-response}
and~\eqref{eq:sens-score-linearization} are evaluated columnwise instead of
forming a dense Jacobian inverse. Although this sparse implementation keeps memory
usage within the 64\,GB allocation, sensitivity scoring still requires repeated
linearized response calculations and NLP re-solves with fixed device settings.

Tables~\ref{tab:benchmark-pm3} and~\ref{tab:benchmark-full} report the comparison.
Here $T$ denotes total runtime in seconds for each method, and the best metric is bolded within each case. {\em The runtime results highlight that sensitivity-guided rounding is substantially slower} on medium-size cases, fails on IEEE 300, and times out on the larger PEGASE 9241 and
13659 cases under an 8\,h limit for both device configurations, yielding 10/13 successful solves.
Homotopy continuation method and FP-storm succeed on all 13 cases. The remaining results show a useful tradeoff. {\em While FP-storm is fast and often competitive in voltage deviation metrics, its feasibility-driven MILP projection step generally yields
weaker VAR reduction and little or no cost improvement relative to ACOPF
dispatch. The homotopy heuristic attains the best VAR quality on seven of ten cases in the tap
$\pm3$ setting and eight of ten cases in the tap $\pm16$ setting, with the lowest cost across
all reported cases and better scalability than the sensitivity baseline.}

\begin{table}[t]
\centering
\scriptsize
\setlength{\tabcolsep}{2pt}
\caption{Heuristic benchmark for tap $\pm3$, CB 0--3, and DCOPF $\to$ ACPF input.}
\label{tab:benchmark-pm3}
\begin{tabular}{l l rrr r r}
\toprule
\textbf{Case} & \textbf{Method}
  & $\mathrm{MAE}_{v}$
  & $\mathrm{MAE}_{q}$
  & $\mathrm{MAE}_{p}$
  & $\%\!\Delta c$
  & $T$ (s) \\
\midrule
\multirow{3}{*}{118}
  & Hom. VVO    & 0.036 & 37.1 & 5.8 & $\mathbf{-0.05}$ & 8.8 \\
  & FP-storm    & 0.035 & \textbf{36.9} & 6.0 & 0.00 & \textbf{4.9} \\
  & Sensitivity & \textbf{0.034} & 40.1 & \textbf{5.6} & $-$0.02 & 9.3 \\
\midrule
\multirow{3}{*}{1354}
  & Hom. VVO    & 0.087 & \textbf{39.4} & \textbf{13.6} & $\mathbf{-0.21}$ & 37.1 \\
  & FP-storm    & \textbf{0.071} & 71.5 & 15.2 & 0.00 & \textbf{17.2} \\
  & Sensitivity & 0.078 & 67.9 & 14.9 & $-$0.10 & 1426.9 \\
\midrule
\multirow{3}{*}{1888}
  & Hom. VVO    & 0.072 & \textbf{24.9} & \textbf{5.1} & $\mathbf{-2.19}$ & 45.5 \\
  & FP-storm    & \textbf{0.039} & 28.7 & 31.8 & 0.00 & \textbf{32.1} \\
  & Sensitivity & 0.061 & 25.1 & \textbf{5.1} & $-$2.15 & 642.7 \\
\midrule
\multirow{3}{*}{2848}
  & Hom. VVO    & 0.062 & 19.7 & \textbf{3.0} & $\mathbf{-0.06}$ & 62.3 \\
  & FP-storm    & \textbf{0.047} & 20.2 & 3.2 & 0.00 & \textbf{20.0} \\
  & Sensitivity & 0.060 & \textbf{19.5} & \textbf{3.0} & $\mathbf{-0.06}$ & 1153.6 \\
\midrule
\multirow{3}{*}{2869}
  & Hom. VVO    & 0.083 & \textbf{37.0} & 20.4 & $\mathbf{-0.19}$ & 100.1 \\
  & FP-storm    & \textbf{0.070} & 68.6 & \textbf{20.3} & 0.00 & \textbf{39.5} \\
  & Sensitivity & 0.078 & 54.0 & 20.5 & $-$0.15 & 10590.5 \\
\midrule
\multirow{3}{*}{6468}
  & Hom. VVO    & 0.061 & \textbf{28.0} & \textbf{10.9} & $\mathbf{-0.14}$ & \textbf{201.3} \\
  & FP-storm    & \textbf{0.048} & 30.0 & 11.2 & 0.00 & 367.0 \\
  & Sensitivity & 0.058 & 29.1 & \textbf{10.9} & $-$0.11 & 13739.1 \\
\midrule
\multirow{3}{*}{6470}
  & Hom. VVO    & 0.058 & \textbf{22.2} & 7.1 & $\mathbf{-0.15}$ & 229.4 \\
  & FP-storm    & \textbf{0.048} & 24.3 & \textbf{7.0} & 0.00 & \textbf{158.3} \\
  & Sensitivity & 0.056 & 23.9 & 7.1 & $-$0.13 & 10541.1 \\
\midrule
\multirow{3}{*}{6495}
  & Hom. VVO    & 0.059 & \textbf{31.0} & \textbf{8.8} & $\mathbf{-0.29}$ & 230.0 \\
  & FP-storm    & \textbf{0.050} & 33.5 & 9.0 & 0.00 & \textbf{91.5} \\
  & Sensitivity & 0.057 & 33.0 & \textbf{8.8} & $-$0.22 & 9912.9 \\
\midrule
\multirow{3}{*}{6515}
  & Hom. VVO    & 0.059 & \textbf{32.6} & \textbf{9.4} & $\mathbf{-0.20}$ & 247.0 \\
  & FP-storm    & \textbf{0.047} & 37.8 & 9.9 & 0.00 & \textbf{148.8} \\
  & Sensitivity & 0.054 & 36.7 & 9.5 & $-$0.13 & 10596.7 \\
\midrule
\multirow{3}{*}{8387}
  & Hom. VVO    & 0.083 & 80.3 & \textbf{38.9} & $\mathbf{-1.09}$ & \textbf{755.3} \\
  & FP-storm    & \textbf{0.079} & 84.0 & 39.7 & 0.11 & 876.0 \\
  & Sensitivity & 0.083 & \textbf{80.2} & \textbf{38.9} & $\mathbf{-1.09}$ & 13935.4 \\
\bottomrule
\end{tabular}
\end{table}

\begin{table}[t]
\centering
\scriptsize
\setlength{\tabcolsep}{2pt}
\caption{Heuristic benchmark for tap $\pm16$, CB 0--3, and DCOPF $\to$ ACPF input.}
\label{tab:benchmark-full}
\begin{tabular}{l l rrr r r}
\toprule
\textbf{Case} & \textbf{Method}
  & $\mathrm{MAE}_{v}$
  & $\mathrm{MAE}_{q}$
  & $\mathrm{MAE}_{p}$
  & $\%\!\Delta c$
  & $T$ (s) \\
\midrule
\multirow{3}{*}{118}
  & Hom. VVO    & 0.036 & \textbf{36.3} & 5.9 & $\mathbf{-0.04}$ & 9.1 \\
  & FP-storm    & \textbf{0.035} & 36.9 & 6.0 & 0.00 & \textbf{4.4} \\
  & Sensitivity & \textbf{0.035} & 39.3 & \textbf{5.6} & $-$0.03 & 9.5 \\
\midrule
\multirow{3}{*}{1354}
  & Hom. VVO    & 0.087 & \textbf{43.2} & \textbf{14.0} & $\mathbf{-0.24}$ & 39.4 \\
  & FP-storm    & \textbf{0.071} & 71.5 & 15.2 & 0.00 & \textbf{17.6} \\
  & Sensitivity & 0.080 & 69.9 & 14.8 & $-$0.14 & 1517.0 \\
\midrule
\multirow{3}{*}{1888}
  & Hom. VVO    & 0.072 & \textbf{24.0} & \textbf{5.0} & $\mathbf{-2.21}$ & 48.6 \\
  & FP-storm    & \textbf{0.039} & 28.7 & 31.8 & 0.00 & \textbf{32.1} \\
  & Sensitivity & 0.072 & 25.4 & \textbf{5.0} & $-$2.20 & 304.8 \\
\midrule
\multirow{3}{*}{2848}
  & Hom. VVO    & 0.065 & 19.4 & \textbf{3.0} & $\mathbf{-0.08}$ & 65.7 \\
  & FP-storm    & \textbf{0.047} & 20.2 & 3.2 & 0.00 & \textbf{20.1} \\
  & Sensitivity & 0.064 & \textbf{18.9} & \textbf{3.0} & $\mathbf{-0.08}$ & 1169.5 \\
\midrule
\multirow{3}{*}{2869}
  & Hom. VVO    & 0.083 & \textbf{38.3} & 20.8 & $\mathbf{-0.23}$ & 95.4 \\
  & FP-storm    & \textbf{0.070} & 68.6 & \textbf{20.3} & 0.00 & \textbf{39.9} \\
  & Sensitivity & 0.078 & 54.6 & 20.8 & $-$0.19 & 10516.7 \\
\midrule
\multirow{3}{*}{6468}
  & Hom. VVO    & 0.062 & \textbf{26.0} & \textbf{10.9} & $\mathbf{-0.19}$ & \textbf{213.4} \\
  & FP-storm    & \textbf{0.048} & 30.0 & 11.2 & 0.00 & 360.2 \\
  & Sensitivity & 0.059 & 28.3 & \textbf{10.9} & $-$0.15 & 16325.3 \\
\midrule
\multirow{3}{*}{6470}
  & Hom. VVO    & 0.058 & \textbf{21.8} & 7.2 & $\mathbf{-0.23}$ & 194.7 \\
  & FP-storm    & \textbf{0.048} & 24.3 & \textbf{7.0} & 0.00 & \textbf{141.6} \\
  & Sensitivity & 0.056 & 23.2 & 7.2 & $-$0.22 & 8207.7 \\
\midrule
\multirow{3}{*}{6495}
  & Hom. VVO    & 0.058 & \textbf{30.5} & 9.1 & $\mathbf{-0.57}$ & 229.1 \\
  & FP-storm    & \textbf{0.050} & 33.5 & \textbf{9.0} & 0.00 & \textbf{82.8} \\
  & Sensitivity & 0.053 & 34.2 & \textbf{9.0} & $-$0.31 & 13529.2 \\
\midrule
\multirow{3}{*}{6515}
  & Hom. VVO    & 0.060 & \textbf{29.8} & 9.6 & $\mathbf{-0.32}$ & 228.1 \\
  & FP-storm    & \textbf{0.047} & 37.8 & 9.9 & 0.00 & \textbf{131.0} \\
  & Sensitivity & 0.054 & 36.1 & \textbf{9.5} & $-$0.16 & 13663.5 \\
\midrule
\multirow{3}{*}{8387}
  & Hom. VVO    & \textbf{0.078} & 81.9 & \textbf{39.6} & $\mathbf{-2.71}$ & \textbf{834.9} \\
  & FP-storm    & 0.079 & 84.0 & 39.7 & 0.11 & 899.8 \\
  & Sensitivity & \textbf{0.078} & \textbf{81.8} & \textbf{39.6} & $-$2.69 & 13190.1 \\
\bottomrule
\end{tabular}
\end{table}

\section{Conclusion}\label{sec:conclusion}

This work shows that transmission-level VVO can serve as a practical
corrective layer between market dispatch and grid operation with improved
Volt/VAR control and lower operational cost.
More broadly, it shows that strategic coordination of existing control devices can support more reliable and efficient operation without requiring large-scale
infrastructure investments.

The weighted objective allows operators to select tradeoffs that reflect their
operational priorities, which is important in practice because system operators
face different ownership structures, reliability requirements, and economic
incentives.
The homotopy continuation method makes these tradeoffs implementable by producing
AC-feasible setpoints across all test cases while maintaining runtime
stability and scalability.

Future work will explore extensions to security-constrained formulations, 
multi-period planning, and learning-augmented frameworks that accelerate optimization
through warm starts or proxy models.

\section*{Acknowledgments}\vspace{-4pt}
This work was supported by the National Science Foundation under Grant
No.~2112533. Any opinions, findings, and conclusions or recommendations
expressed in this material are those of the author(s) and do not necessarily
reflect the views of the sponsors.

\printbibliography

\end{document}